\documentclass{article}
\usepackage{arxiv}
\usepackage[utf8]{inputenc} 
\usepackage[T1]{fontenc}    
\usepackage{hyperref}       
\usepackage{url}            
\usepackage{booktabs}       
\usepackage{amsfonts}       
\usepackage{nicefrac}       
\usepackage{microtype}      

\usepackage{latexsym}
\usepackage{amsmath}
\usepackage{amssymb}
\usepackage[all]{xy}
\usepackage[mathscr]{eucal}

\newcommand{\ptes}{\mathcal P}

\newcommand{\caf}{\mathcal F}

\usepackage{graphicx}
\graphicspath{ {./images/} }

\newcommand{\esp}{\vspace*{0.5cm}}
\newcommand{\esple}{\vspace*{0.3cm}}
\newcommand{\vaz}{\emptyset}
\newcommand{\s}{\subseteq}

\newcommand{\co}{\supseteq}

\newcommand{\meni}{\leqslant}

\newcommand{\mai}{\geqslant}

\newcommand{\pe}{\langle}
\newcommand{\pd}{\rangle}
\newcommand{\dominado}{\preccurlyeq}

\newcommand{\pv}{\mathbf{\mathcal P \mathcal V}}
\newcommand{\dialc}{\textrm{Dial}(\textrm{\bf Sets})}
\newcommand{\dial}{\textrm{Dial}_2(\textrm{\bf Sets})^{\textrm{op}}}
\newcommand{\dialo}{\textrm{Dial}_2(\textrm{\bf Sets})}
\newcommand{\mhd}{\mathbf{MHD}}
\newcommand{\zf}{\mathbf{ZF}}
\newcommand{\zfc}{\mathbf{ZFC}}
\newcommand{\ac}{\mathbf{AC}}

\newcommand{\acw}{\mathbf{AC}_\omega}
\newcommand{\dc}{\mathbf{DC}}
\newcommand{\menorgt}{\meni_{GT}}

\newcommand{\pvunb}{\pv_{\textrm{unb}}}

\newcommand{\R}{\mathbb{R}}
\newcommand{\N}{\mathbb{N}}

\newtheorem{Df}{Definition}[section]
\newtheorem{Th}[Df]{Theorem}
\newtheorem{Prop}[Df]{Proposition}

\newtheorem{Ex}[Df]{Example}

\newtheorem{Que}[Df]{Question}

\newcommand{\fd}{\hfill{$\blacksquare$}}
\newcommand{\bc}{\begin{center}}
\newcommand{\ec}{\end{center}}
\newcommand{\n}{\noindent}
\newcommand{\dem}{{\bf Proof:}}

\title{Kolmogorov-Veloso Problems\\
and Dialectica Categories}

\author{
Valeria de Paiva \\
Topos Institute\\
Berkeley, CA \\
\texttt{valeria@topos.institute} \\
\And
Samuel G. da Silva \\
Departamento de Matem\'{a}tica \\
Instituto de Matem\'{a}tica e Estat\'{i}stica\\
Universidade Federal da Bahia\\
\texttt{samuel@ufba.br}
}

\begin{document}
\maketitle

\begin{abstract}
We investigate the categorical connection between Dialectica constructions, Kolmogorov problems, Veloso problems and Blass problems. We  show that the work of Kolmogorov can be regarded as a bridge between  Veloso abstract notion of a problem and the conceptual problems Blass discussed in his questions-and-answers framework. This bridge can be seen by means of the categorical Dialectica constructions introduced in de Paiva's dissertation and reformulated by da Silva to account for set-theoretical foundational assumptions. The use of categorical concepts  allows us to provide several examples, connecting extremely different areas of mathematics, and using simple methods. This paper also shows that while Blass’ and Kolmogorov’s notions of problem can be investigated using  Zermelo-Fraenkel's (ZF) set-theoretical framework, Veloso’s problems require the Axiom of Choice (AC). Moreover, weaker notions of choice (dependent choice and countable choice) can also be accounted for in the problems' framework.
\end{abstract}

\section{Introduction}
Blass' seminal paper \cite{blass1995} establishes a surprising connection between de Paiva's Dialectica categories \cite{dePaiva1991}, Vojt\'a\v s' methods to prove inequalities between cardinal characteristics of the continuum~\cite{vojtas1993}  and the complexity theoretical notions of
problems and reductions developed in~\cite{levin1}. 
Blass does not mention Kolmogorov's very abstract notion of \textit{problem} (\cite{kolmogorov1932}), which is not related to specific complexity issues. Kolmogorov did investigate  a notion of abstract problem, producing  an alternative intuitive semantics for Propositional Intuitionistic Logic, part of the celebrated $BHK$ interpretation of Intuitionistic Logic.
 
 Kolmogorov's problems were cited as an inspiration by Veloso when he developed his own theory of abstract problems in the Eighties \cite{veloso1984}, but the two frameworks  were not formally connected. 
 This note recalls  Veloso's `Teoria de Problemas' (Theory of Problems) and shows how it can be related to the Dialectica construction~\cite{dePaiva1989}, via  Kolmogorov's concepts.  To establish this relationship, we use the  modification of the Dialectica construction considered by the second author \cite{silva2017, silva2020} for other  set-theoretical purposes. This modification corresponds to a condition of non-triviality of the collections of problems or solutions, first suggested by \cite{moore2004}.

The categorical connection between Dialectica models, Kolmogorov's problems, Veloso's problems and Blass' problems shows that the use of categories really allows us to connect extremely different areas of mathematics, using simple methods. In the case of this paper, it also allows us to determine where exactly the foundational choices of axioms are important. We show that while Blass' and Kolmogorov's notions of problem can be investigated using the set-theoretical framework of $\zf$, Veloso's problems commit us to a stronger set-theory, as discussed in the following sections. Whether this requirement of stronger foundations is a bug or a feature depends, perhaps, on personal taste and conviction. From our part we are happy to note that, like many other questions in Mathematics, as soon as one investigates them a little, the notion of `problem' points out to grand challenges in the foundations of Mathematical Logic, to wit whether one wants to accept or not the Axiom of Choice within their chosen framework.

\section{Kolmogorov Problems}
Blass (\cite{blass1995}) noticed that de Paiva's Dialectica construction $GC$ \cite{dePaiva1989}, when the base category $C$ is the category $\bf Sets$, is the dual of Vojt\'a\v s’s category $GT$ of generalized Galois-Tukey connections
\cite{vojtas1993}. Since the Dialectica construction has been generalized in many different directions, instead of writing $G{\bf Sets}$ we will write $\dialo$ for this category, to make explicit the object $2$ (the object of truth-values) where our relations  map into, as well as the category $C$ that is ${\bf Sets}$, where  objects `live'. We recall the definition of the category below.

\begin{Df} [Dialectica category \cite{dePaiva1989}] The category $\dialo$ has as objects triples of the form $A = (U,X,\alpha)$, where $U$ and $X$ are sets and $\alpha \s U \times X$ is a set-theoretical relation, which can be written, equivalently, as $\alpha: U\times X\to 2$. 
If $A = (U,X,\alpha)$ and $B = (V,Y,\beta)$ are objects of $\dialo$, a morphism from $A$ to $B$ is a pair of functions in ${\bf Sets}$, $(f,F)$, $f:U \to V$, $F:Y \to X$ such that
\bc For all $u \in U$ and $y \in Y$, $u \alpha F(y)$ implies $f(u)\beta y$. \ec \end{Df}

 Blass also noticed  that 
the opposite, dual Dialectica Category $\dial$ can be taken to intuitively mean that objects represent \textit{problems} and morphisms stand for reductions between problems. 
Thus a triple $P = (I,S,\sigma)$ of $\dial$, under this interpretation, represents a problem, whose instances are elements of $I$, the set of possible solutions is given by $S$ and the relation $\sigma$ can be read as  ``is solved by", that is, if $z$ is an instance of the problem $P$ and $s$ is a possible solution then $z \sigma s$ states that ``$s$ solves $z$". 
Blass associates these problems with a concept of many-one reduction of search problems in complexity theory, a very restricted kind of problem, for  which he refers to \cite{levin1, levin2}.

Veloso's theory of problems (\cite{veloso1984}), following P\'olya, suggests 
that in order to understand a problem one should consider the following  initial questions:
\begin{enumerate}
\item What is the unknown? 
\item What are the data? 
\item What is the condition?
\end{enumerate}
These questions correspond directly to the elements of the triples of information which characterize a problem, which will be referred to, in this work, as a {\it Kolgomorov problem}. They  also correspond precisely to Blass' interpretation above.

\esple

\begin{Df}[Kolgomorov problems] A {\bf Kolgomorov problem} is a triple $P = (I, S, \sigma)$, where $I$ and $S$ are sets and $\sigma \s I \times S$ is a  set theoretical relation. We  say that:
\begin{itemize}
\item $I$ is the set of {\bf instances} of the problem $P$;

\item $S$ is the set of {\bf possible solutions} for the instances $I$; and

\item $\sigma$ is the {\bf problem condition}, i.e. the relation $\sigma$ holds between $z$ and $s$, in symbols $z \,\sigma\, s$ 
if  the solution $s$ satisfies the problem condition $\sigma$ for the instance $z$, or, more briefly,  ``$s$ $\sigma$-solves $z$". 
\end{itemize}
\end{Df}

\esple

Kolgomorov's 1932 paper, republished in English in 1991 (\cite{kolmogorov1932}) has two parts. The first section, which introduces Kolgomorov's problems, says that ``If the intuitionistic cognitive presuppositions are not accepted then one should take into account only the first section".
In this section he  introduces problems via mathematical examples, for instance ``Find  any four integers, $x,y,z$ and $n$ such that $x^n+ y^n= z^n$, for $n > 2$". (Note that Kolmogorov states explicitly, in page 152, that ``We never assume a problem to be solvable".)

We can identify the {\it objects} of the category $\dialo$ with  Kolgomorov problems. What would the {\it morphisms} represent in this case? To answer this, we
consider the morphisms of the opposite of the Dialectica category, $\dial$. In $\dial$,  a morphism from an object $P' = (I', S', \sigma')$ to an object $P =(I, S,\sigma)$ is a pair of functions $(f, F)$, where  $f\colon  I \to I'$ and 
$F\colon S'  \to S $ are such that the following condition holds
\esple

\begin{center}
$(\forall z\in I) \ (\forall t\in S') \ [f(z)\,\sigma'\,t \longrightarrow z \,\sigma\, F(t)].$
\end{center} 
 So for all instances of problems $z$ of $P$ and all solutions of problems $t$ of $P'$, if  the instance of problem $f(z)$ has $\sigma'$-solution $t$ then the instance $z$ has $\sigma$-solution $F(t)$.
 
 If we regard $P$ and $P'$ as Kolgomorov problems, the existence of a morphism from $P'$ to $P$ ensures that there is a {\it reduction} of the problem $P$ to the problem $P'$ -- because the act of solving an instance of $P$ may be reduced to the act of solving an instance of $P'$. More precisely, if one wants to solve a particular instance $z$ of the Kolgomorov problem $P$, it suffices to find a solution $t$ for the instance $f(z)$ of $P'$ -- since $F(t)$ will provide a solution for the initial instance $z$ of $P$. 
 
 Kolmogorov discusses a few number-theoretical and geometrical problems, as well as abstract, logical rules  ones. In what follows, we present a number of examples from daily mathematical practice to show how they are coded as Kolgomorov problems. We also show reductions between those problems which can be seen as morphisms of $\dial$. As a piece of notation, if $X$ is a set and $n \in \N$ then $[X]^n$ denotes the family of all subsets of $X$ which have precisely $n$ elements and  $[X]^{\meni n}$  means $\bigcup\limits_{m \meni n} [X]^m$.
 \esple
 
 \begin{Ex} Analytical geometry is entirely based on the reduction of geometrical problems to equation solving problems.
 \end{Ex}
 
We present some practical examples to explain what we mean by the statement above. 

Let $\pi$ be any plane of the 3-dimensional Euclidean space $\R^3$ and let $F\colon\R^2 \to \pi$ be a coordinate system as usual (i.e., for every pair $(u,v)$ of real numbers we associate the point $P = F(u,v)$ of the plane which has coordinates $(u,v)$ -- that is, $P = P_{(u,v)}$). Every line $l$ of the plane $\pi$ is then represented by an equation of the form $ax + by = c$, where $a$ and $b$ are real numbers such that $a \neq 0$ or $b \neq 0$ and $c \in \{0,1\}$. Let $E$ be the family of all equations of the described canonical form and let $\mathcal{L}$ denote the family of all lines of the plane $\pi$. The decision problem of ``whether a given point lies on a given line" is $(\mathcal{L},\pi,\ni)$, and the problem of ``whether a given pair of real numbers satisfy a given equation" is $(E,\R^2,\zeta)$, on which an equation $ax + by = c$ is $\zeta$-related to a pair $(u,v)$ of real numbers if $au + bv = c$. Then we can reduce the geometrical problem $(\mathcal{L},\pi,\ni)$ to the algebraic problem $(E,\R^2,\zeta)$ using the morphism $(f,F)$, where $f: \mathcal{L} \to E$ is defined by putting  $f(l) = \textrm{eq}(l)$ (where $\textrm{eq}(l)$ is the canonical equation which represents $l$) and $F: \R^2 \to \pi$ is the coordinate system. Indeed, if $(u,v)$ satisfies the equation of a line $l$ we know that the corresponding point $P_{(u,v)}$ lies in $l$.
 
 A slight variation of what we have just done reduces the problem of finding the intersection point of two distinct lines to the problem of solving a linear system with two equations over two variables. The geometrical problem of finding the intersection point of two distinct lines is $([\mathcal{L}]^2,\pi,\xi)$, where $\{l_1,l_2\}\,\xi\,P$ means that $P \in l_1 \cap l_2$ for any distinct lines $l_1, l_2$ of $\pi$ and for every point $P$ in $\pi$. The algebraic problem of finding the solution for a linear system with two equations over two variables is $([E]^2,\R^2,\lambda)$, where the relation $\lambda$ in $\{a_1x + b_1y = c_1, a_2x + b_2y = c_2\} \lambda (u,v)$ is the conjunction of $``a_1x + b_1y = c_1" \zeta (u,v)$ and $``a_2x + b_2y = c_2" \zeta (u,v)$. The morphism which gives the reduction is $(g,F)$, where $g:[\mathcal{L}]^2 \to [E]^2$ is given by $g(\{l_1,l_2\}) = \{f(l_1),f(l_2)\} = \{\textrm{eq}(l_1),\textrm{eq}(l_2)\}$ for all distinct lines $l_1,l_2$ of $\pi$ and $F:\R^2 \to \pi$ is still the coordinate system. Now, if $l_1$ and $l_2$ are distinct lines and $(u,v)$ is a solution of the system $\{\textrm{eq}(l_1),\textrm{eq}(l_2)\}$ then the point $P_{(u,v)}$ lies in the intersection of the lines $l_1$ and $l_2$.
 
 \esp
 
 \begin{Ex} The problem of finding vectors in the intersection of the kernels of a finite family of linear functionals over $\R^n$ reduces to the problem of finding the solutions of a homogeneous system of linear equations. 
 \end{Ex}
 
 Let $n \mai 2$ and $\{e_1,e_2,\ldots,e_n\}$ denote the canonical basis of the $n$-dimensional Euclidean space $\R^n$, regarded as a vector space over $\R$. A {\it linear functional} over $\R^n$ is a linear transformation from $\R^n$ into $\R$.  Let
 $(\R^n)^*$ (the {\it dual} of $\R^n$) denote the family of all linear functionals over $\R^n$. It is well-known that the dimension of the dual space is also $n$, and that $\{\mathbf{x_1},\mathbf{x_2},\ldots,\mathbf{x_n}\}$ is a basis of the dual space -- where, for $1 \meni i \meni n$,  $\mathbf{x_i}: \R^n \to \R$ is the linear functional which satisfies $\mathbf{x_i}(e_j) = \delta_{ij}$ (where $\delta_{ij} = 1$ if $i = j$ and it is zero otherwise) for $1 \meni j \meni n$  and is then extended to all linear functionals by linearity. 
 
 Let $\mathcal{E}$ denote the family of all linear equations of the form $a_1x_1 + a_2x_2 + \ldots + a_nx_n = 0$, where $a_1,\ldots,a_n$ are real numbers. A  linear functional   $\mathbf{T}$ in the dual space of $\R^n$ may be written, in a unique way, in the form $\mathbf{T} = \sum\limits_{i = 1}^{n} a_i\mathbf{x_i}$.

 Thus we may define a translation function $g:(\R^n)^* \to \mathcal{E}$ in the obvious way, i.e. by putting $$g(\mathbf{T}) = ``a_1x_1 + a_2x_2 + \ldots + a_nx_n = 0"$$ if  $\mathbf{T} = \sum\limits_{i = 1}^{n} a_i\mathbf{x_i}$.
 
 The problem of finding vectors in the intersection of the kernel of finite families of linear functionals is $([(\R^n)^*]^{\meni n},\R^n,\alpha)$, where $\mathcal{H}\,\alpha\,(c_1,c_2,\ldots,c_n)$ means that $(c_1,c_2,\ldots,c_n) \in \bigcap\limits_{h \in \mathcal{H}} \textrm{Ker}(h)$ for any finite family $\mathcal{H}$ of $m$ linear functionals with $m \meni n$ and any $n$-tuple $(c_1,c_2,\ldots,c_n) \in \R^n$ (recall that the dimension of the dual space is $n$ as well, so we may only consider finite families with no more than $n$ linear functionals). The problem of solving a homogeneous system of no more than $n$ linear equations is $([\mathcal{E}]^{\meni n},\beta,\R^n)$, where $\mathcal{G}\,\beta\,(c_1,c_2,\ldots,c_n)$ means that, for any $(c_1,c_2,\ldots,c_2)\in \R^n$ and any $\mathcal{G}\in[E]^m$ with $m \meni n$,  
 $(c_1,c_2,\ldots,c_2)$ solves each one of the $m$ linear equations of the finite family $\mathcal{G}$. A morphism which codes the usual and easily verifiable procedure of reducing the first problem to the second is now given by $(f,Id)$, where $f:[(\R^n)^*]^{\meni n} \to [\mathcal{E}]^{\meni n}$ is given by $f(\mathcal{H}) = \{g(h): h \in \mathcal{H}\}$ for any finite family $\mathcal{H}$ of no more than $n$ linear functionals and $Id$ is the identity function from the {\it set}\, $\R^n$ into the {\it vector space}\, $\R^n$ -- i.e., we have reduced a problem stated in the context of the linear structure of $\R^n$ and its dual to a more naive, pedestrian problem of finding  solutions of a homogeneous system of equations.

 Veloso developed further notions of {\it viable problems}, {\it links}  between problems (and {\it reduction links} between problems) in his theory of mathematical problems \cite{veloso1984}. However, as Veloso himself was aware of (see further discussion in the final section),  these definitions make essential use of the Axiom of Choice (\cite{martin-loef2006},\cite{moore1982}), and this use will be detailed and discussed in the next section. 

\section{Veloso Problems and the Axiom of Choice}
In our previous work~\cite{silva2017} we described a modification of the 
category $\dial$, which Blass calls the category $\cal{PV}$
in \cite{blass1995}. 
Blass, describing the category $\cal{PV}$, explains that it has as objects problems, together with their instances and respective solutions. Moreover, morphisms of the category $\cal{PV}$  identify reductions of classes of (complexity) problems to others. The modification of the category in \cite{silva2017}, following the work of \cite{moore2004}, insists on some conditions of non-triviality of the objects.
These conditions can be paraphrased as `there are no problems without a solution' and `there are no solutions that solve all problems at once'. 
There are also some constraints  on the cardinality of the sets of instances and solutions of $\pv$; all constituent sets of objects of $\pv$ are bounded above by the cardinality of the {\it continuum} (which is $\mathfrak{c} = 2^{\aleph_0} = |\R| = |\ptes(\N)|$).

The mentioned non-triviality conditions are not found in the original Dialectica categorical construction (\cite{dePaiva1989},\cite{dePaiva1991}), where trivial objects are actually required to provide truth-values (or units for the categorical operators) associated to the logical connectives of Girard's Linear Logic~\cite{girard1987}. The non-triviality conditions are also not found in Kolmogorov's work. He describes \textit{meaningless} problems as the ones that do not have a solution.

The non-triviality conditions
will characterize, exactly, a class of problems we refer to as {\it Veloso problems}. We note, however, that there are no upper bounds for the cardinality of the sets of instances and solutions of the problems originally discussed by Veloso.

We  first formally present the category $\pv$, as well as stratified versions of it which were introduced in \cite{silva2020}.
The following definition should be considered within $\zfc$, since we refer to well-ordered cardinals in its first clause.  

\begin{Df} [{Category $\pv$ \cite{silva2017}}]\label{pv} The category $\pv$ is the subcategory of $\dial$ whose objects are the triples $A = (U,X,\alpha)$ satisfying the following three clauses, which will be referred to as the  $\mhd$ conditions ($\mhd$ stands for Moore, Hru\v s\'ak and D\v zamonja~\cite{moore2004}):

\esple

(1) The cardinalities of the (non-empty) constituent sets  are bounded above by the cardinality of the continuum --  i.e. $0 < |U|,|X| \meni 2^{\aleph_0}$.  

\esple

(2) Every problem has a solution -- i.e. $$(\forall u \in U)(\exists x \in X)[u\,\alpha\, x].$$

\esple

(3) There are no solutions that solve all the problems at once -- i.e. $$(\forall x \in X)(\exists u \in U)[\neg\, (u\, \alpha\, x)].$$

\esple

\normalsize

The morphisms between objects of $\pv$ are the  same morphisms of $\dial$ -- that is, a morphism from an object $B = (V, Y, \beta)$ to an object $A =(U,X,\alpha)$ is a pair of functions $(f, F)$, where  $f:  U \to V$ and 
$F: Y \to X$ are such that 
\esple

\begin{center}
$(\forall u\in U) \ (\forall y\in Y) \ [f(u)\,\beta\,y \longrightarrow u \,\alpha\, F(y)].$
\end{center} \end{Df}
\esple
Morphisms of $\pv$ induce the 
{\bf Galois-Tukey pre-order} introduced by Vojt\'a\v s, which is defined in the following way:  if $A = (U,X,\alpha)$ and $B = (V,Y,\beta)$ are objects of $\mathcal{PV}$,  then we have

\bc $A \meni_{GT} B  \Longleftrightarrow \mbox{  There is a morphism from } B \mbox{ to } A$. \ec
\esple

The diagram below represents the situation where  $A \meni_{GT} B$:

\[
\xymatrix{
&^u \, U \ar[d]_{f}& \alpha &  \,\,\,\,\,\,\,\,X^{F(y)}\\
& _{f(u)}\, V\,\,\,\,\,\, & \beta & \,\,\,{Y}\, _{y} \ar[u]_{F}}
\]

Given an object $A = (U,X,\alpha)$ of $\pv$, its {\it dual object} is given by $A^* = (X,U,\alpha^*)$, where $x \alpha^* u$ means that $\neg\,(u \alpha x)$. One can easily check (via a contrapositive argument) that: \bc If $A \menorgt B$, then $B^* \menorgt A^*$. \ec

In \cite{silva2020} parametrized, stratified versions of $\pv$ were introduced -- the $\pv_X$ categories, for $X$ an infinite set.
The main goal of the proposed stratification was to generalize features of $\pv$ to sets of higher cardinalities, since the constituents of the objects of $\pv$ are bounded above by the cardinality of the  continuum (which is the cardinality of the power set of the naturals). Notice that, indeed, in $\zfc$ the categories $\pv$ and  $\pv_{\N}$ coincide, and so the whole idea of the stratified versions was to generalize $\pv$ (in a choiceless context) 
to any other other infinite set $X$.
These categories are introduced in the $\zf$ setting, as their definitions use the notion of \textit{domination} ($\dominado $) between sets instead of that of well-ordered cardinals (for a more detailed discussion on the subtleties of comparing sizes in the absence of $\ac$, we refer to Section 2 of \cite{silva2020}). Recall that a set $A$ is {\it dominated} by a set $B$, $A \dominado B$, if there is an injective function from $A$ into $B$.

\esple

\begin{Df}[{Categories $\pv_X$}]\label{pvX} Let $X$ be an infinite set in $\zf$. Then $\pv_X$ is the subcategory of $\dial$ whose objects $(U,X,\alpha)$ are those which satisfy the following $\mhd_X$ conditions (since
they are a restricted form of the $\mhd$ conditions):

\esple

(1) $U$, $V$ are non-empty sets and $U,V \dominado \ptes(X)$, that is, sets $U,V$ are dominated by $\ptes(X)$. 

\esple

(2) $(\forall u \in U)(\exists x \in X)[u\,\alpha\, x]$

\esple

(3) $(\forall x \in X)(\exists u \in U)[\neg\,(u\,\alpha\,x)]$

\esple

The morphisms between objects of $\pv_X$ are the same morphisms of $\dial$. 

\end{Df}

\esple

We also define, in the expected way, a Galois-Tukey ordering $\meni_{GT}$ on the objects of $\pv_X$.

The categories $\pv_X$ are used in \cite{silva2020} to provide, 
after a quantification over all infinite sets, equivalences of the Axiom of Choice $\ac$. The following equivalences are proved in \cite{silva2020}:

\begin{Th}[da Silva {\cite{silva2020}}] The Axiom of Choice is equivalent to the following statements:

\esple

\n (*) {\it  ``For every infinite set $X$, $(\ptes(X),\ptes(X),=)$ is a maximum element in the Galois-Tukey ordering $\meni_{GT}$ on  $\pv_X$."\,}

\esple

\n and

\esple

\n (**) {\it  ``For every infinite set $X$, $(\ptes(X),\ptes(X),\neq)$ is a minimum element in the Galois-Tukey ordering $\meni_{GT}$ on  $\pv_X$."\,}
\end{Th}

\esple

This theorem shows that stratifying objects and features of the Dialectica construction using the Galois-Tukey ordering is equivalent to accepting the Axiom of Choice. This may be a surprise to mathematicians not used to  ``thinking about foundational issues".

In this work we are interested in the investigation of general problems, which means that, in particular, the sets of instances and solutions should have no upper bounds on their cardinalities. To do so in a categorical setting, we provide the following definition, which
should be considered in the  setting of $\zf$.

\esple

\begin{Df}[{Unbounded $\pv$ category}] \label{pvunb} The category $\pvunb$, the {\it unbounded $\pv$ category}, is the the subcategory of $\dial$ 
whose objects $(U,X,\alpha)$ satisfy the $\mhd$ conditions (2) and (3) but where the $\mhd$ condition (1) is relaxed to allow sets of arbitrary large cardinality. Thus we require $\mhd$ conditions (2) and (3) together with

\esple

(1)'  $U$ and $X$ are non-empty sets

\esple

A triple $(U,X,\alpha)$ is an object of $\pvunb$ if it is an object of $\pv_Y$ for some infinite set $Y$.

\esple

The morphisms between objects of $\pvunb$ are the very same morphisms of $\dial$. 

\end{Df}

\esple

The unbounded $\pv$ category captures, precisely {\it Veloso problems} -- these are the Kolgomorov problems which are {\bf viable} and {\bf non-generic}, as we define now.

\esple

\begin{Df}[Veloso problem] Let $P = (I,S,\sigma)$ be a Kolgomorov problem. 

$(i)$ $P$ is said to be a {\bf viable} problem if the domain of $\sigma$ is the whole set $I$ -- or, equivalently, if for every instance $z$ of $I$ there is some possible solution $s$ in $S$ that solves the problem, so that the relation $z\,\sigma\,s$ holds.
          
$(ii)$ A possible solution $s$ will be said to be a {\bf generic solution} for $P$ if it solves all of its instances  -- that is, if for every instance $z$ of $I$ one has $z\,\sigma\,s$ for this particular solution $s$.

$(iii)$ The problem $P$ is said to be {\bf non-generic} if it has no generic solutions -- equivalently, for every possible solution $s \in S$ there is some instance $z \in I$ such that $\neg (z \sigma s)$.

$(iv)$ $P$ will be said to be a {\bf Veloso problem} if it is viable and non-generic.

\end{Df}

\esple

Notice that the fact that  a problem $P = (I,S,\sigma)$  is non-generic is equivalent to the viability of the dual problem $P^* = (S,I,\sigma^*)$, where $s\,\sigma^*\,z$ means $\neg (z \sigma s)$ for all $z \in I$ and $s \in S$.

It should be clear that viability and non-genericity are conditions asking for  non-triviality of a given problem. For instance, it is well-known that the following statement is an Axiom of Incidence, in Hilbert's plane geometry:

\esple

($\dagger$) {\it ``For every line $l$, there are at least two points which lie in the line $l$ and at least one point which does not lie in the line $l$"}

\esple

Together with the axiom ``For every pair of distinct points there is only one line which passes through those points", the above statement ($\dagger$) corresponds, precisely, to the viability and non-genericity requirements for the problem of determining whether a given line contains a given pair of distinct points of the plane or $(L,P,\co)$, where 

\esple

$L$ = the set of all lines of the plane; and

\esple

$P$ = the family of all pairs of two distinct points of the plane.

\esple

Notice also that non-genericity establishes a dimension for the preceding problem: if all points were in the same line $l$ we would not be studying the geometry of the plane, only the geometry of a line. 

The following example shows that certain mathematical notions are equivalent to the viability of certain problems. Recall that a subset of a topological space is {\it dense} if its closure is equal to the whole space. It is a textbook easy exercise to show that dense sets are precisely those which intersect any non-empty open set. Thus, the following holds:

\esple

\begin{Ex} Let $(X,\tau)$ be topological space and $D$ be a proper subset of $X$. Then $D$ is a dense subspace of $X$ if, and only if, the Kolgomorov problem $(\tau \setminus \{\vaz\},D,\ni)$ is viable. 
\end{Ex}

\esple

Alternatively, $D$ is a dense set if, and only if, for every non-empty open set $U$ of $X$ the problem $(\{U\},D,\ni)$ is viable. This observation is relevant in the following example, which is, accordingly to Veloso (\cite{veloso1984}, page 24), a {\it problem to prove} in  P\'olya's terminology  -- in contrast to the so-called {\it problems to find}. Recall that a topological space is a {\it Baire space} if the intersection of any countable family of dense subsets of the space is also dense subset of the space.

\esple

\begin{Ex}\label{baireantes} To prove the Baire Theorem for Complete Metric Spaces (i.e., the theorem which asserts that any complete metric space is a Baire space), it suffices to verify that, for any non-empty open set $O$ and for any countable family $\{U_n:n \in \N\}$ of dense open sets, the problem $(\{O\},\bigcap\limits_{n \in \N} U_n,\ni)$ is a viable problem. 
\end{Ex}

\esple

Later on we will show that the above example is easier to handle if transformed into a countable chain of viable related problems. 

The notion of \textit{viability} of a problem, introduced by Veloso in \cite{veloso1984} (page 25), is  related to the notion of {\it solvability} of the problem. However, a solution of a problem for Veloso is represented by a function, not by a relation, as it is the case in our work. We will use the terminology {\it solution function} to denote these objects that are, instead of solution relations, solution functions. 

\esple

\begin{Df} Let $P = (I,S,\sigma)$ be any viable Kolgomorov problem. A {\bf solution function} for $P$ is a function $f:I \to S$ satisfying $f \s \sigma$ -- or, equivalently, $f$ is a $S$-valued function with domain $I$ 
such that for every instance $z$ of $P$ we have that $f(z)$ solves $z$.
\end{Df}

\esple

Veloso (op. cit) has defined a notion of solvability according to his definition of solution -- that is, in our terms,  a Kolgomorov problem $P$ is {\bf solvable} if $P$ has a solution function. As one of the easiest equivalent statements of the Axiom of Choice is precisely {\it ``Every relation contains a function with the same domain"} (see page 131 of \cite{moore1982}),
it is straightforward to check that, within $\zfc$ (i.e., assuming the Axiom of Choice), the following equivalences hold:

\esple

\begin{Prop}[viable is solvable in $\zfc$] \label{viavelsoluvel}  Let $P = (I,S,\sigma)$ be any Kolgomorov problem. The following statements are equivalent:

\esple

$(i)$ The Kolmogorov problem $P$ is viable.

$(ii)$  The problem $P$ is solvable.

$(iii)$ The problem $P$ satisfies the formula $(\forall z \in I)(\exists s \in S)[z\,\sigma\,s]$.
\end{Prop}

The preceding proposition was also stated  by
Veloso.
The only non-obvious implication ($(i) \longrightarrow (ii)$) follows easily from the  equivalence of the Axiom of Choice mentioned above.

Still following Veloso's definitions from \cite{veloso1984}, next we define the notions of {\it links} and {\it reduction links} between problems. We  point out that the notion of ``reduction link"  will be defined in terms of ``solution functions". Again, this context is intrinsically associated to the Axiom of Choice, as we will discuss later.

\esple

\begin{Df}[Links and reduction links \cite{veloso1984}]\label{links}
Let $P = (I,X,\sigma)$ and $P' = (I',X',\sigma')$ be two Kolgomorov problems. 

\esple

$(i)$ A {\bf link} from $P$ to $P'$ is a pair of functions $(f,F)$, where

\esple

$f:I \to I'$ is said to be a {\bf translation function}; and

$F: S' \to S$ is said to be a {\bf recovery function}.

\esple

$(ii)$ A link $(f,F)$ from $P$ to $P'$ is called a {\bf reduction link} of $P$ to $P'$ if it lifts solution functions from $P'$ to $P$, i.e. for every solution function $g$ of $P'$ the composite function $F \circ g \circ f$ is a solution function for $P$.

\end{Df}

\esple

From now on we  identify the class of all Veloso problems with the objects of the category $\pvunb$. Next we show that morphisms of this category correspond to reduction links between the corresponding problems.

\esple

\begin{Prop} \label{morlink} Let $P = (I,X,\sigma)$ and $P' = (I',X',\sigma')$ be objects of $\pvunb$, considered as Kolgomorov problems. If $(f,F)$ is a morphism witnessing the order $P \meni_{GT} P'$, then $(f,F)$ is a reduction link of $P$ to $P'$. 
\end{Prop}

\n \dem. Let $g$ be an arbitrary solution function for the problem $P'$. 
As we know that $(f,F)$ is a morphism from $P'$ to $P$, we know that
for all $z \in I$ and $t \in S'$ the following implication holds:

\bc $f(z)\,\sigma'\,t \longrightarrow z\,\sigma\,F(t)$ \ec

Fix an arbitrary $z \in I$. As $g \s \sigma'$, for $t = g(f(z))$ we have that $f(z)\,\sigma'\,t$, and therefore
$$z \sigma F(t) = F(g(f(z))$$
and thus $(z, F(g(f(z)) \in \sigma$. By the arbitrariness of $z$, we conclude that $F \circ g \circ f \s \sigma$ and so the composite function $F \circ g \circ f$ is a solution function for $P$, as desired. 
\fd

\esple

So, given Kolgomorov problems $P$ and $P'$ with a Dialectica morphism $(f,F)$ witnessing the inequality $P \meni_{GT} P'$ in the Galois-Tukey ordering, we are allowed to interpret such inequality as a measure of complexity, since a solution of $P$ may be reduced to a solution of $P'$ -- thus under $P \meni_{GT} P'$ we can say that $P$ is as easy to solve as $P'$, or that $P$ is not more complicated to be solved than $P'$.

We have shown that the existence of solution functions for viable problems is ensured by the Axiom of Choice (Proposition \ref{viavelsoluvel}).
As Veloso's  reduction links were defined in terms of solution functions, it is clear that this approach presupposes a choice principle. In the following results, we will show  more,  as it will be established that assuming that Veloso's approach holds in full generality gives rise to a number of equivalents of the Axiom of Choice (or of weak related statements).  For the rest of this section, our base theory is $\zf$.

\esple

\begin{Th}\label{teo1} $\ac$ is equivalent to the following statement:

\bc ``Every Veloso problem has a solution function"\,. \ec \end{Th}

\n \dem. ($\Rightarrow$) Let $P = (I,S,\sigma)$ be a Veloso problem. By the $\mhd$ condition (1), $P$ is viable and so $dom(\sigma) = I$. As in the proof of Proposition \ref{viavelsoluvel},  the existence of a solution function is guaranteed by the equivalent of $\ac$ given by the statement ``Every relation contains a function with same domain". 

\n ($\Leftarrow$) It was remarked in \cite{silva2020} that $\ac$ is equivalent to the following statement:

\esple

{\it ``For every infinite set $X$, there is a choice function defined on $\ptes(X) \setminus \{\vaz\}$ -- that is, there is $f: \ptes(X) \setminus \{\vaz\} \to X$ such that $f(Y) \in Y$ for all $\vaz \neq Y \s X$.''}

\esple

So, let $X$ be any infinite set. The problem 
 $$(\ptes(X) \setminus \{\vaz\}, X, \ni)$$ is clearly viable (since non-empty subsets have elements, by definition) and non-generic (given $x \in X$ one has that $Y = X \setminus \{x\}$ is not solved by $x$). A solution function for this problem is, clearly, a choice function 
for  $\ptes(X) \setminus \{\vaz\}$. As $X$ was taken arbitrarily, we have established $\ac$ by the above remark. \fd

\esple

It follows from the previous theorem that in the absence of the Axiom of Choice there will be Veloso problems in $\pv_{\textrm{unb}}$ without any solution function. 

Notice that, if  $P'$ is a Veloso problem without a solution function, it is vacuously true that  {\it any link} from  $P$ to $P'$ (for {\it any} given $P$) is a reduction link. Motivated by this, we introduce the following definition:

\esple

\begin{Df} Let $P, P'$ be problems in $\pv_{\textrm{unb}}$. A morphism $(f,F)$ from $P'$ to $P$ is said to be {\bf realized as a reduction link} if the statement ``$(f,F)$ is a reduction link from $P$ to P'\,\,"\, holds non-vacuously. 

\end{Df}

If  all morphisms have to be realized as reduction links, then the Axiom of Choice must be present, as the following theorem shows. 

\begin{Th} \label{teo2}
 $\ac$ is equivalent to the statement

\bc ``Every morphism of $\pv_{\textrm{unb}}$\,\,is realized as a reduction link\,".  \ec \end{Th}

\n \dem. ($\Rightarrow$) Given a morphism from $P'$ to $P$, $\ac$ implies that $P'$ has solution functions (by \ref{teo1}), and the rest follows from Proposition \ref{morlink} and from the definition of ``realized as a reduction link". 

\n ($\Leftarrow$) Given any infinite set $X$, we are able to consider the object of $\pvunb$ given by
$$P = (\ptes(X) \setminus \{\vaz\}, X, \ni).$$ and the corresponding identity morphism 
$(id,id)$ from $P$ to $P$. If we assume that $(id,id)$ is realized as a reduction link then there is a solution link for $P$, which will be a choice function for $$\ptes(X) \setminus \{\vaz\}.$$ As $X$ was taken arbitrarily, we obtain $\ac$ in the same way as  in Theorem \ref{teo1} \fd

\esple

There are several set-theoretical statements which are referred to, in the literature, as {\it weak choice principles}. Weak choice principles are implied by the Axiom of Choice $\ac$, but they are not equivalent to it, and are often regarded as ``fragments" or ``partial cases" of $\ac$. To identify the precise amount of choice which is needed for a particular argument/result is a fruitful and current line of research (akin to Reverse Mathematics) within Set Theory. This may be seen in the standard reference \cite{HowardRubin} and in the dozens of papers which have cited it over the last twenty years. 

The Axiom of Countable Choice (usually denoted by $\acw$) is one of the most celebrated weak choice principles.  It corresponds to the restriction of the Axiom of Choice to countable families of non-empty sets, that is, it asserts that ``Every countable family of non-empty sets has a choice function", or if $\{X_n: n \in \N\}$ if a countable family of non-empty sets then there is a function $$f:\{X_n: n \in \N\} \to \bigcup\limits_{n \in \N} X_n $$ such that $f(X_n) \in X_n$ for all $n \in \N$. 

In the following theorem, we  prove an equivalent of countable choice $\acw$ in terms of Veloso problems -- more specifically, in terms of Veloso problems whose instance sets are countable.

\begin{Th}\label{teo3}  The Axiom of Countable Choice $\acw$ is equivalent to the following statement:

\bc ``Every Veloso problem whose set of instances is countable \\ has a solution function". \ec \end{Th}

\esple

\n \dem. ($\Rightarrow$) Let $P = (I,S,\sigma)$ be a Veloso problem, on which the instance set $I$ is a countable set. Enumerate $I = \{z_n: n \in \N\}$. As $P$ is viable, we know that $$(\forall n \in \N)(\exists s \in S)[z_n\,\,\sigma\,\,s]$$ \n and so for every $n \in \N$ the set $$F_n = \{s \in S: z_n\,\sigma\,s\}$$ \n is a non-empty set. Applying $\mathbf{AC}_\omega$ to the countable family of non-empty sets given by $$\caf = \{F_n: n \in \N\},$$ \n we get a choice function $g: \caf \to \bigcup \caf$ with $g \s \sigma$ -- so that $z_n\,\,\sigma\,\,g(F_n)$ for all $n \in \N$. A solution function $f$ for the problem $P = (I,S,\sigma)$ is now easily defined by putting, for every $z \in I$, $$f(z) = g(F_m) \iff z = z_m.$$

\n ($\Leftarrow$) Let $\caf = \{X_n: n \in \N\}$ be a countable family of non-empty sets. We have to show that there is a choice function for such family, i.e. we have to exhibit a function $f: \caf \to \bigcup \caf$ such that $f(X_n) \in X_n$ for all $n \in \N$.

Consider the set $X$ given by $$X = \bigcup\limits_{n \in \N}\, (\{n\} \times X_n) $$ \n and let $\sigma \s \N \times X$ be the binary relation defined in the following way: 

$$n\,\sigma\, x \iff \Pi_1(x) = n,$$

\n where $\Pi_1$ denotes the projection on the first coordinate. In other words, for all $(m,z) \in  X$ we have that $n\,\sigma\,(m,z)$ if, and only if, $n = m$.

Consider the problem $(\N,X,\sigma)$. As $\caf$ is supposed to be a family of non-empty sets, such problem is viable, and (as $\N$ is infinite) it is easy to check that it is also non-generic. So, $(\N,X,\sigma)$ is a Veloso problem. By hypothesis, there is a solution function of $(\N,X,\sigma)$, say $g: \N \to X$. Now we define $$f:\caf \to \bigcup\limits_{n \in \N} X_n$$ \n by putting $$f(F_n) = \Pi_2(g(n))$$ \n for all $n \in \N$. It should be clear that $f$ is a choice function for $\caf$. \fd

\esple

The reader may check that we could have proved versions of the three preceding theorems stated only in terms of viable problems. 
In fact, Veloso himself has observed that the equivalence between $\ac$ and the statement ``Every viable problem has a solution" holds (\cite{veloso1984}, page 35). However, as the problem $(\ptes(X) \setminus \{\vaz\},X, \ni)$ is a Veloso problem (i.e., viable and non-generic) for any infinite set $X$ -- and it is, basically, the problem which appears in several parts of the proofs -, we have preferred to point out that the class of Veloso problems is enough, in each case, to provide an equivalence with the corresponding choice principle.

The final result of this section will be stated in terms of viable problems only. First, we will introduce the notion of {\it $\omega$-chain of viable problems}. 

\esple

\begin{Df}[$\omega$-chain of viable problems] Let $\{I_n: n \in \N\}$ be a family of non-empty sets and $\{\sigma_n: n \in \N\}$ be a family of binary relations such that, for all $n \in \N$,  $$\sigma_n \s I_n \times I_{n + 1}.$$ Let $S_n = I_{n + 1}$ and $P_n = (I_n,S_n,\sigma_n)$ for all $n \in \N$. We say that $\pe P_n: n \in \N\pd$ is an {\bf $\omega$-chain of viable problems} if all problems $P_n$ are viable, i.e. for all $n \in \N$ one has $$\forall z \in I_n\,\, \exists s \in S_n\,\, [z\,\sigma_n\, s].$$ \end{Df} 

\esple

The preceding definition aims to capture the very common situation in Mathematics where one has to solve a $\N$-sequence of chained problems, and the solution of the $n$-th problem is an instance of the $(n+1)$-problem. In such situation one wishes, of course, to produce a sequence of solutions. 

\esple 

\begin{Df} Let $\pe P_n: n \in \N\pd$ be an $\omega$-chain of viable problems, where 
$P_n = (I_n,S_n,\sigma_n)$ and $S_n = I_{n+1}$ for all $n \in \N$. A {\bf solution sequence} for $\pe P_n: n \in \N \pd$ is a sequence $\pe z_n: n \in \N\pd$ such that $z_n \in I_n$ and  $z_n \sigma_n z_{n + 1}$ for all $n \in \N$.
\end{Df}

\esple

We present below a mathematical example where the above concepts play an important role. 

\begin{Ex} The usual proof of  the Baire Category Theorem for complete metric spaces given in the textbooks may be encoded by a $\omega$-chain of viable problems. \end{Ex}

Let us see why. Let $\{U_n: n \in \N\}$ be a countable family of dense open sets of a complete metric space $M$ and let $O$ be an arbitrary non-empty open set (as in Example \ref{baireantes}). We have to show that the intersection of all the $U_n$'s meets $O$. We let $I_n$ be the family of all non-empty open sets whose diameter is less than $\frac{1}{n+1}$ and whose closures are included in $U_n \cap O$, and let $\sigma_n$ be the reverse inclusion. It is easy to check that all of the problems $\pe I_n,S_n,\sigma_n\pd$ (where $S_n = I_{n + 1}$) are viable problems. If $\pe V_n: n \in \N\pd$ is a solution sequence for such $\omega$-chain of problems, then $\big(\bigcap\limits_{n \in \N} U_n) \cap 0$ is ensured to be non-empty, since $\{\overline{V_n}: n \in \N\}$ is a decreasing sequence of non-empty closed sets whose diameters converge to zero, so it has non-empty intersection in the complete metric space $M$ --  and, by the definition of the $I_n$'s, any point of such intersection testifies that $\big(\bigcap\limits_{n \in \N} U_n) \cap 0 \neq \vaz$.

What is usually not said in the textbooks about the preceding example is that a solution sequence for the $\omega$-chain of viable problems is given by the {\it Principle of Dependent Choices}, denoted by $\dc$, which is another celebrated weak choice principle\footnote{In fact, the Principle of Dependent Choices is {\it equivalent} to the Baire Theorem for complete metric spaces, as shown in \cite{blair}.}. The Principle of Dependent Choices states:  

\esple
 {\it ``If $\delta$ is a binary relation on a non-empty set $A$ (i.e. $\delta \s A \times A$) satisfying $$(\forall x \in A)(\exists y \in A)[x\,\delta\,y],$$ then there is a sequence $\pe x_n :n \in \N\pd$ of elements of $A$ such that $x_n\,\delta\,x_{n + 1}$ for all $n \in \N$''.}
 
 \esple

The principle $\mathbf{DC}$ is regarded as the precise amount of choice needed to make a countable number of consecutive arbitrary choices. It is well-known that the Axiom of Choice is stronger than the Principle of Dependent Choice which is stronger than the Axiom of Countable Choice, i.e.
$$\mathbf{AC} \Rightarrow \mathbf{DC} \Rightarrow \acw,$$ and that none of those implications is reversible. More information on these and many other weak choice principles may be found in \cite{HowardRubin}. 

In what follows, we present equivalents of $\dc$ in terms of $\omega$-chains of viable problems. 

\begin{Th} \label{teo4} The following statements are equivalent:

\esple

$(i)$ The Principle of Dependent Choices, $\dc$;

$(ii)$ For any $\omega$-chain of viable problems  $\pe P_n: n \in \N\pd$ (where, for all $n \in \N$,
$P_n = (I_n,S_n,\sigma_n)$ and $S_n = I_{n+1}$) and for every instance $z \in I_0$, there is a solution sequence $\pe z_n:n \in \N \pd$ for this $\omega$-chain with $z_0 = z$; and

$(iii)$ Every $\omega$-chain of viable problems has a solution sequence. 

\esple

\end{Th}

\n \dem.$(i) \Rightarrow (ii)$: Let $\pe P_n: n \in \N\pd$ be as in the statement. We define a set $A$ whose elements are, precisely, all $(k+1)$-tuples $$(z_0,z_1,\ldots,z_k),$$ where $k$ ranges over all natural numbers, $z_0 = z$, $z_i \in I_i$ for all $0 \meni i \meni k$ and also $z_i\,\sigma_i\,z_{i+1}$ for all $0 \meni i \meni k - 1$ if $k > 0$. We define, over the tuples of $A$, a relation $\delta$ such that $$(w_0,w_1,\ldots,w_k)\,\delta\,(t_0,t_1,\ldots,t_j)$$ if $k < j$ and $w_i = t_i$ for all $0 \meni i \meni k$ -- that is, $\delta$ is the usual strict prefix order over the tuples. 

As all of the problems of the $\omega$-chain are viable by definition, it is easy to check that the set $A$ and the relation $\delta$ satisfy the requirements one needs to apply the Principle of Dependent Choices. So, by $\dc$, there is a sequence $\pe s_n: n \in \N\pd$ of tuples such that $s_n\,\delta\,s_{n+1}$ for all $n \in \N$. If we identify each tuple with the corresponding finite sequence, we get a compatible family of functions (in fact, an increasing chain of compatible finite functions), and so the union of such family is a function. It is easy to check that the obtained sequence $z = \bigcup\limits_{n \in \N} s_n$ is a solution sequence for the $\omega$-chain with $z_0 = z$, as desired.

\n $(ii) \Rightarrow (iii)$: Obvious.

\n $(iii) \Rightarrow (i)$: If $A$ and $\delta$ are under the hypothesis of $\dc$, then $A$ is a non-empty set and $(A,A,\delta)$ is a viable problem. Then we just have to consider the $\omega$-chain $\pe P_n: n \in \N\pd$ where $I_n = S_n = A$ and $\sigma_n = \delta$ for all $n \in \N$. By $(iii)$ we may assume that there is a solution sequence for this $\omega$-chain, say $\pe z_n: n \in \N\pd$. So, this sequence satisfies $z_n\,\delta\,z_{n + 1}$ for all $n \in \N$. As $A$ and $\delta$ were taken arbitrarily, we have just established $\dc$ -- and the proof is then finished. \fd

\esple

We have shown that the Axiom of Countable Choice, the Principle of Dependent Choices and the fully-fledged Axiom of Choice correspond to natural (sub-)classes of problems, as described by Kolmogorov and Veloso.

\section{Tight Coupled Problems}
In this short section we discuss how the notion of a ``tight coupled reduction link" between problems can be seen as a previously known variant of the Dialectica construction we have discussed so far.

The way in which the notion of {\it reduction link} of $P$ to $P'$ was defined (item $(ii)$ of Definition \ref{links}) seems to suggest that there is no tight connection or ``coupling", in general, between problems $P$ and $P'$ -- as pointed out by Veloso in the page 28 of \cite{veloso1984}. Veloso argues that a reduction link of $P$ to $P'$ is, in most cases, {\it uncoupled} in the following sense: after applying the translation function $f$ on an instance $z$ of $P$, we may forget completely about $P$ and only care about solving the instance $f(z)$ of $P'$ -- and,  after that, any solution $t \in S'$ of the instance $f(z)$ of $P'$ will generate a solution of the instance $z$ of $P$ by applying the recovery function $F$. However, 
in certain situations (see example below), it is interesting to allow the use of some additional information during the recovery process. 

For this case, Veloso 
defines  {\it tightly coupled links} from $P$ to $P'$ in the following way: the translation function is still  some $f:I \to I'$ but the recovery function is a function  $F\colon I \times S' \to S$ so that the information on the original instance $z$ of $P$ may be used at the time of recovery -- that is, the pair $(f,F)$ needs to satisfy the following condition: for any $z \in I$ and any $t \in S'$, $$f(z)\,\sigma'\,t\longrightarrow   z\,\sigma\,F(z,t).$$ 

This  requirement on the  tightly coupled links corresponds, precisely, to the morphisms of the dual of $\dialc$ (the simplest case of the ``original" Dialectica category, inspired by G\"odel's Dialectica Interpretation and introduced by the first author in \cite{dePaiva1989:AMS}).

\begin{Ex}[{Reduction with tight coupled link}] Using tight coupled links, one can formally prove that the problem of finding a normal line of a surface through a certain point is not more complicated than the problem of finding orthogonal planes to a certain plane. 
\end{Ex}

In this example, we use {\it surface} for a two-dimensional differential manifold $S \s \R^3$, and therefore given a point $x \in S$ there is a tangent plane $T_x(S)$ centered at $x$. Recall that if $x$ is a point of a plane $\pi \s \R^3$ then a line $l$ through $x$ is the {\it normal line} of the plane $\pi$ (through $x$) if $l$ is perpendicular to all lines of $\pi$ which go through $x$, and if $x$ is a point of a surface $S$ then the normal line of $S$ through $x$ is the normal line of the tangent plane $T_x(S)$ through $x$. Two intersecting planes $\pi$ and $\rho$ are said to be {\it orthogonal} if for every point of the intersection line the normal lines of $\pi$ are contained in $\rho$ and vice-versa. Notice that if $l$ is the normal line of a plane $\pi$ through $x$ then $l$ is included in every plane $\rho$ which is orthogonal to $\pi$ and passes through $x$ -- in other words, all planes which are orthogonal to a given plane $\pi$ and passes through a given point share the normal line of $\pi$ through this very same point.

Let $\mathcal{L}$ be the family of all lines of the 3-dimensional Euclidean space $\R^3$ and $\mathcal{P}$ the family of all planes of $\R^3$. For a given surface $S$, the problem of finding the normal lines through each point of the surface is $(S,\mathcal{L},\sigma)$, where $x\,\sigma\,l$ means ``$l$ is the normal line of $S$ through $x$" for every point $x \in S$ and every line  $l \in \mathcal{L}$, and the problem of finding orthogonal planes is $(\mathcal{P},\mathcal{P},\xi)$, where $\pi\, \xi\, \rho$ means ``$\pi$ and $\rho$ are orthogonal planes" for all planes $\pi$ and $\rho$ in $\R^3$. 

A tight coupled link from $(S,\mathcal{L},\sigma)$ to $(\mathcal{P},\mathcal{P},\xi)$ is given by the pair $(f,F)$, where $f: S \to \mathcal{P}$ is defined by putting $f(x) = T_x(S)$ for all $x \in S$ and $F\colon S \times \mathcal{P} \to \mathcal{L}$ is given by $F(x,\rho) = \varphi(T_x(S),\rho,x)$, where $\varphi: \mathcal{P} \times \mathcal{P} \times \R^3 \to \mathcal{L}$ is defined by putting $\varphi(\pi,\rho,x) =$ the line $l$ contained in the plane $t(x,\rho)$ (where $t(x,\rho)$ is $\rho$ itself if $x \in \rho$ or is the unique plane parallel to $\rho$ passing through $x$ otherwise) which is perpendicular to the intersecting line of $\pi$ and $t(x,\rho)$ through $x$, if $\pi$ and $\rho$ are orthogonal planes; and any previously fixed line otherwise. In view of the fact remarked at the end of the previous paragraph, it is easy to see that $F(x,\rho)$ is the normal line of $T_x(S)$ through $x$ (and thus, of $S$) whenever $\rho$ is a plane orthogonal to $f(x) = T_x(S)$. Notice that the described reduction link formalizes the following mental procedure: ``if I know how to produce an orthogonal plane for any given plane, then I know how to produce the normal line of a surface $S$ at a point $x$: we take any plane $\rho$ that is orthogonal to $T_x(S)$, translate it to $x$ -- via parallel translation -- and then consider the perpendicular line (through $x$ and contained in the translated plane) of the intersection line of $T_x(S)$ and the translated plane".

This shows one case where the information on the original problem instance (the original point $x$) is required for the recovery function $F$ of the tight coupled link $(f,F)$ that reduces the original problem of ``finding normal lines to a given point in the surface" to the new problem of ``finding orthogonal planes". This problem is particularly nice, as it seems to connect to approaches to automatic differentiation under development using categorical machinery in \cite{elliott2018}. More research work is required here.

\section{Conclusions and Further Work}

We have shown that the work of Kolmogorov can be regarded as a bridge between the abstract problems Veloso discussed in his Theory of Mathematical Problems and the complexity problems Blass discussed in \cite{blass1995}. This bridge can be seen by means of the categorical Dialectica constructions $GC$ and $DC$ introduced in \cite{dePaiva1991}.

 Veloso's {\it Critical Retrospect} in \cite{veloso1984}  already observed that his approach on viable Kolgomorov problems  relies heavily on the Axiom of Choice 
 and asks whether such approach (together with some suggested techniques on decomposition of problems) was ``tainted"\, by $\ac$ -- since the mathematical entities whose existence depend on the Axiom of Choice (and statements about those entities) are the usual examples of non-constructive notions in Mathematics, and he aimed  his approach to incorporate constructivity at some level. In this paper we have proceed with this line of research and we introduced a restricted class of Kolgomorov problems, which we call {\it Veloso Problems}, and we have shown that if we restrict ourselves to this class and assume that the presented machinery holds in full generality then we get, again, equivalences of the Axiom of Choice. We have also shown that some related/restricted notions on problems are intrinsically associated to weak choice principles such as the Axiom of Countable Choice and the Principle of Dependent Choices. To determine precisely the deductive strength of assertions relating to Kolgomorov and Veloso problems (positioning them in the hierarchy of weak choice principles) seems to deserve further research, and the same could be said about the following question: 

\begin{Que}[{implicit in \cite{veloso1984}}] How to deal with the above mentioned non-constructive aspects of Kolgomorov-Veloso theory of problems within constructive environments? \end{Que}

The work in this paper shows that some answers to this  question are obtained by focusing on {\it relations}, instead of {\it functions}. In fact, equivalences with choice principles arise from assuming  either the existence of solution functions (Theorems \ref{teo1}, \ref{teo2} and \ref{teo3}) or of solution sequences (Theorem \ref{teo4}). 

Veloso (op.cit.) presents some arguments against the use of relations instead of functions in the representation of solutions of problems. 
He argues that if the solution of a problem can be represented as a relation, then the problem condition itself, $\sigma$, would be a solution of the problem, and so nothing more would need to be done; to know the specification of the problem would solve it automatically. 
Second, he argues that, if one wants to assume that some instance $z \in I$ could be solved by more than one element of the set of possible solutions $S$ 
then, even in this case, it is possible to work with a representation using functions, but in this case one should work with the so-called {\it multifunctions}, or multivalued functions. If $I$ is the set of instances and $S$ is the set of possible solutions of a viable problem $P$, a multifunction solution $f$ would be, formally, a function with domain $I$ and codomain $\ptes(S)$ (i.e, the set of all subsets of $S$) such that, for every $z \in I$, $f(z)$ is the {\it subset } of $S$ given by $\{s \in S: z\,\sigma\,s\}$ -- that is, the use of multivalued solution functions consists in associating each instance $z \in I$ to the set of {\it all} of its particular solutions. 

Replying to the first argument, in practical applications, to know the {\it definition} of the condition problem $\sigma$ does not give a solution automatically (even less the set of all solutions) for an instance $z \in I$. In simple/naive cases, the more realistic approach would be, probably, to proceed with some decision problem/search problem considering, for each $z \in I$,
the set $\{z\}\times S$ as a domain and  with $(\{z\}\times S) \cap \sigma$ as the set of yes-instances for this problem. So, to know the problem condition does not solve the problem. For instance, one may consider the viable problem $(\{\zeta\},\mathbb{C},\sigma)$, on which $\zeta$ is Riemann's zeta function and $\sigma$ is the relation given by $\zeta\,\sigma\, c \iff \zeta(c) = 0$ -- that is, $\sigma$ is the restriction to $\{\zeta\}$ of the very general problem of finding zeros of analytic functions. Despite  $\sigma$ being perfectly (and easily) defined, we still do not know (after more than 150 years, see \cite{Bombieri}) whether there are complex numbers (apart from the even negative numbers) with real part distinct from $\frac{1}{2}$ which solve the problem. 

For the second objection against relations as solutions, we would like to mention two  points. The first is that, if we decide to work with the powerset $\ptes(S)$ instead of $S$, then we are dramatically increasing the cardinality of the sets we are dealing with. For problems where any instance has only a {\it finite} set of solutions there is perhaps no major issue
but we want to work with  theories where  problems have instances with possible infinite solutions. (Actually the set-theory applications do insist on infinite sets.) As some features of the theory may depend on the cardinalities of the constituents of the triples of the corresponding categories (as in the definitions of  $\pv$ and $\pv_X$ for a given infinite set $X$, see Definitions 
\ref{pv} and \ref{pvX}), it would not be desirable such a dramatic increase in cardinality. The second argument is that if we assume $f(z)$ to be the set of all solutions (for $f$ a multivalued solution function and $z$ some instance of the problem) then we have all solutions for $z$ ``locked inside a box"; these solutions may become individually inaccessible and this precludes a qualitative analysis of them. We prefer to have the possibility of comparing distinct solutions for any fixed instance of the problem we are solving. Thus we insist that relations are a better modelling tool than functions, one that allows us to move on to functions, when and if we feel $\ac$ is adequate.

A second main conclusion for us is that the notions of reduction between problems, considered by Kolmogorov, Veloso and Blass, seen to be well modelled by the categorical morphisms in either $\dialo$, its dual or the original Dialectica construction. A skeptical reader may complain that the categorical language used is not buying us much. We beg to differ: the possibility of relating formally these, to begin with, quite `woolly' notions of problems and solutions, seems a serious step forward in the hard task of detecting unwarranted foundational assumptions that tend to `sneak' into mathematics. We still need to investigate whether the traditional tools of category theory, e.g. products, coproducts, exponentials, (co-)limits, etc. can be leveraged to our advantage. And, apart from such traditional tools, we are also interested in the investigation of the possible interactions between the Dialectica categorical modelling of problems and the so-called {\it lenses} (\cite{hedges_2016}, \cite{spivak_2019}). Lenses are constructions used in situations where some structure is converted to different forms -- through actions and observations between environments and agents -- in such a way that all changes made can be reflected as updates to the original structure. Such constructions
 have  attracted the attention of several researchers over the past ten years (see also \cite{nlab:lens_(in_computer_science)}, and references therein). 

A possible avenue for further work from this point on would be to connect the categorical semantics meaning of $\dialo$ (logically speaking $\dialo$ models Linear Logic, together with Intuitionistic Propositional Logic) to, yet to be conceived, models of Ecumenical Propositional Logic. Ecumenical Propositional Logic \cite{pimentel2019} is Prawitz' recent suggestion of how to consider under the same umbrella both intuitionsitic and classical principles, as used by mathematicians. Since both Kolmogorov and Veloso mentioned their intentions of being understood by both classical and intuitionistic mathematicians, it would be extremely nice if the categorical models discussed here could help with modelling ecumenical logic. However, new insights will be required to deal with the traditional issues of modelling categorically classical logic.

\bibliographystyle{unsrt}  
\bibliography{references1}  

\end{document}